\documentclass[a4paper]{ifacconf}

\usepackage{graphicx}      
\usepackage{natbib}        

\usepackage{makeidx,epsfig}
\usepackage{amsmath}
\usepackage{amsfonts}
\usepackage{amssymb}
\usepackage{multirow}
\usepackage{mathrsfs}

\newtheorem{corollary}{Corollary}
\newtheorem{definition}{Definition}
\newtheorem{example}{Example}
\newtheorem{proposition}{Proposition}
\newtheorem{theorem}{Theorem}
\begin{document}
\begin{frontmatter}

\title{Fault Tolerant Stabilizability\\ of Multi-Hop Control Networks\thanksref{footnoteinfo}}

\thanks[footnoteinfo]{This research has been partially supported by the European Commission through the HYCON2 project.}

\author{M.D. Di Benedetto},
\author{A. D'Innocenzo}, {\bf and}
\author{E. Serra}

\address{Department of Electrical and Information Engineering,\\
Center of Excellence DEWS. University of L'Aquila - L'Aquila, Italy.\\
\small{email:}\{mariadomenica.dibenedetto,alessandro.dinnocenzo,emmanuele.serra\}@univaq.it}

\begin{abstract}
A Multi-hop Control Network (MCN) consists of a plant where the communication between sensor, actuator and computational unit is supported by a wireless multi-hop communication network, and data flow is performed using scheduling and routing of sensing and actuation data. We address the problem of characterizing controllability and observability of a MCN, by means of necessary and sufficient conditions on the plant dynamics and on the communication scheduling and routing. We provide a methodology to design scheduling and routing, in order to satisfy controllability and observability of a MCN for any fault occurrence in a given set of failures configurations.
\end{abstract}


\end{frontmatter}


\baselineskip10.3pt

\section{Introduction} \label{secIntro}

Wireless networked control systems are spatially distributed control systems where the communication between sensors, actuators, and computational units is supported by a shared wireless communication network. Control with wireless technologies typically involves multiple communication hops for conveying information from sensors to the controller and from the controller to actuators. The use of wireless networked control systems in industrial automation results in flexible architectures and generally reduces installation, debugging, diagnostic and maintenance costs with respect to wired networks. The main motivation for studying such systems is the emerging use of wireless technologies in control systems (see e.g.,~\cite{akyildiz_wireless_2004}, \cite{song_wirelesshart:_2008}, \cite{song_complete_2008}).

Although multi-hop networks offer many advantages, their use for control is a challenge when one has to take into account the joint dynamics of the plant and of the communication protocol. Wide deployment of wireless industrial automation requires substantial progress in wireless transmission, networking and control, in order to provide formal models and verification/design methodologies for wireless networked control system. The design of the control system has to take into account the presence of the network, as it represents the interconnection between the plant and the controller, and thus affects the dynamical behavior of the system. Using a wireless communication medium, new issues such as fading and time-varying throughput in communication channels have to be addressed, and communication delays and packet losses may occur. Moreover analysis of stability, performance, and reliability of real implementations of wireless networked control systems requires addressing issues such as scheduling and routing using real communication protocols.

While most of the research on networked control systems is on direct networking, we focus on multi-hop networks. In Section \ref{secRelatedWork} we relate our research to the existing scientific literature on Networked Control Systems. In particular, the modeling and stability verification problem for a MIMO LTI plant embedded in a multi-hop control network (MCN) when the controller is already designed has been addressed in \cite{AlurRTAS09}. A mathematical framework has been proposed, that allows modeling the MAC layer (communication scheduling) and the Network layer (routing) of the recently developed wireless industrial control protocols, such as WirelessHART (\texttt{www.hartcomm2.org}) and ISA-100 (\texttt{www.isa.org}). The mathematical framework defined in \cite{AlurRTAS09} is compositional, namely it is possible to exploit compositional operators of automata to design scalable scheduling and routing for multiple control loops closed on the same multi-hop control network.

In this paper, starting from the mathematical framework developed in \cite{AlurRTAS09}, we address the novel issue of characterizing controllability and observability of a continuous-time SISO LTI plant embedded in a MCN that implements scheduling and routing protocols, and where failures of communication links may occur. We motivate the exploitation of redundancy in data communication (i.e. sending sensing and actuation data through multiple paths) with the aim of rendering the system robust with respect to link failures (e.g. when the battery of a node discharges or a communication channel goes down), and to mitigate the effect of packet losses (e.g. transmission errors).

In Section \ref{secMCNModeling} we extend the model in \cite{AlurRTAS09} to model redundancy, by defining a weight function that specifies how the duplicate information transmitted through the multi-hop network is merged, and by defining a semantics of the redundant data flow through the network. Of course, all results stated in this paper also apply to MCN that do not exploit redundancy. We remark that the differences introduced in this paper with respect to the model in \cite{AlurRTAS09} do not invalidate compositionality of the framework.

As a first result of this paper, given a MCN, we state in Section~\ref{secMCNStabilizability} necessary and sufficient controllability and observability conditions on the plant dynamics and on the scheduling and routing of the communication network.

As a second result, given a MCN and a set of failures configurations of the communication nodes, we state in Section~\ref{secMCNRobustStabilizability} necessary and sufficient conditions on the plant dynamics, on the scheduling and routing of the communication network, and on the set of failures configurations, such that there exists a scheduling and routing configuration that guarantees reachability and observability conditions of the MCN for each failures configuration. Since we adopt a constructive proof, we provide a methodology to configure scheduling and routing of a MCN, in order to satisfy controllability and observability of the closed loop system for any fault occurrence in a given set of failures configurations.

In Section \ref{secSimulations} we apply our results, combined with fault detection and hybrid observer techniques, in order to perform co-design of control algorithms and communication parameters for stabilizing MCN where failures of links occur.

\section{Related work} \label{secRelatedWork}

There exists a wide literature on Networked Control Systems, see for example \cite{Zhang2001}, \cite{WalshCSM2001}, \cite{SpecialIssueNCS2004}, \cite{Hespanha2007} and references therein. The literature on robust stability of networked control systems (see e.g. \cite{hai_lin_robust_2003}, \cite{cloosterman_robust_2006}, \cite{shi_towards_2006}) generally addresses stability analysis in presence of packet loss and variable delays, but does not take into account the non--idealities introduced by scheduling and routing communication protocols of multi-hop control networks. When relating our paper with the current research about the interaction of control networks and communication protocols, most efforts in the literature focus on scheduling message and sampling time assignment for sensors/actuators and controllers interconnected by wired common-bus networks, e.g. \cite{Astrom97j1}, \cite{walsh_stability_2002}, \cite{yook_trading_2002}, \cite{TabbaraCDC2007}, \cite{TabbaraTAC2007}. The authors in \cite{WitrantMSC2007} use model predictive control to stabilize a plant over a multi-hop control network, by only considering delay introduced by the routing policy.

However, what is needed for modeling and analyzing control protocols on multi hop control networks is an integrated framework for analysing/co-designing network topology, scheduling, routing, transmission errors and control. To the best of our knowledge, the only formal model of multi-hop wireless sensor and actuator networks is reported in~\cite{andersson_simulation_2005}. In this paper, a simulation environment that facilitates simulation of computer nodes and communication networks interacting with the continuous-time dynamics of the real world is presented. The main difference between the work presented in~\cite{andersson_simulation_2005} and this work is that here we provide results on a formal mathematical model that takes into account plant dynamics and scheduling-routing dynamics.

At the best of our knowledge, our work is pioneering in addressing the controller design problem for multi-hop control networks that implement standardized scheduling and routing communication protocols, in order to enable co-design of controller, scheduling and routing.


\section{Modeling of MCNs} \label{secMCNModeling}

\begin{figure}[ht]
\begin{center}
\includegraphics[width=0.5\textwidth]{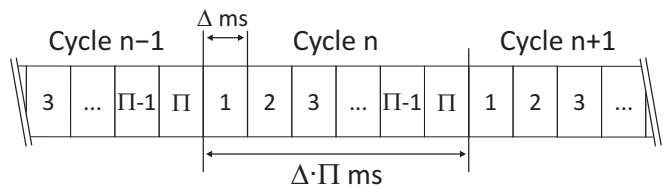}
\caption{Time-slotted structure of the scheduling period.}\label{frame}
\end{center}
\end{figure}

The challenges in modeling multi-hop control networks are best explained by considering the recently developed wireless industrial control protocols, such as WirelessHART and ISA-100. These standards allow designers of wireless control networks to distribute a synchronous communication schedule to all communication nodes of a wireless network. For each working frequency, time is divided into slots of fixed time length $\Delta$ (see Figure \ref{frame}). A periodic scheduling composed by $\Pi$ time slots allows each node to transmit data only in a subset of time slots and frequencies, i.e. a mixed TDMA and FDMA MAC protocol is used. The standard specifies a syntax for defining schedules and a mechanism to apply them. However, the issue of designing schedules and routing remains a challenge for the engineers, and is currently done using heuristic rules. To allow systematic methods for designing schedules that preserve controllability and observability of a plant, a mathematical model of the effect of scheduling and routing on the control system is needed. The MCN model we propose in this paper allows modeling multi-hop control networks that implement the protocols WirelessHART and ISA-100, but it is much more general: it allows modeling general routing and scheduling communication protocols that specify TDMA, FDMA and/or CDMA access to a shared communication resource, for a set of communication nodes interconnected by an arbitrary radio connectivity graph.
\begin{definition}\label{defMCN}
A SISO Multi-hop control network (MCN) is a tuple $\mathcal N = ( \mathcal P, \mathcal G_{\mathcal R}, \eta_{\mathcal R}, \mathcal G_{\mathcal O}, \eta_{\mathcal O}, \Delta)$ where:
\begin{itemize}
\item $\mathcal P = (A,B,C)$ models a plant dynamics in terms of matrices of a continuous-time SISO LTI system.

\item $\mathcal G_{\mathcal R} = ( V_{\mathcal R},E_{\mathcal R}, W_{\mathcal R} )$ is the controllability radio connectivity acyclic graph, where the vertices correspond to the nodes of the network, and an edge from $v_1$ to $v_2$ means that $v_2$ can receive messages transmitted by $v_1$ through the wireless communication link $(v_1,v_2)$. We denote $v_{c}$ the special node of $V_{\mathcal R}$ that corresponds to the controller, and $v_u \in V_{\mathcal R}$ the special node that corresponds to the actuator of the input $u$ of $\mathcal P$. The weight function $W_{\mathcal R} : E_{\mathcal R} \to \mathbb R^+$\footnote{We denote by $\mathbb R^+$ the set of strictly positive reals} associates to each link a positive constant. The semantics of $W_{\mathcal R}$ will be clear in the following definition of $\eta_{\mathcal R}$.

\item $\eta_{\mathcal R} \colon \{1, \ldots, \Pi\} \to 2^{E_{\mathcal R}}$ is the controllability communication scheduling function, that associates to each time slot $k \in \{1, \ldots, \Pi\}$ a set of edges of the controllability radio connectivity graph. The integer constant $\Pi$ is the period of the reachability communication scheduling. The semantics of $\eta_{\mathcal R}$ is that $( v_1,v_2 ) \in \eta(k)$ if and only if at time slot $k$ the data content of the node $v_1$ is transmitted to the node $v_2$, multiplied for the weight $W_{\mathcal R}(v_1,v_2)$.

\item $\mathcal G_{\mathcal O} = ( V_{\mathcal O},E_{\mathcal O}, W_{\mathcal O} )$ is the observability radio connectivity acyclic graph, and is defined similarly to $\mathcal G_{\mathcal R}$. We denote with $v_{c}$\footnote{With abuse of notation, we denote $v_{c}$ both for $\mathcal G_{\mathcal R}$ and $\mathcal G_{\mathcal O}$: this will not lead to confusion, since it will be always clear from the context whether we will be considering the controller node of $\mathcal G_{\mathcal R}$ or $\mathcal G_{\mathcal O}$} the special node of $V_{\mathcal O}$ that corresponds to the controller, and $v_y \in V_{\mathcal O}$ the special node that corresponds to the sensor of the output $y$ of $\mathcal P$.

\item $\eta_{\mathcal O} \colon \{1, \ldots, \Pi\} \to 2^{E_{\mathcal O}}$ is the observability communication scheduling function, and is defined similarly to $\eta_{\mathcal R}$. We remark that $\Pi$ is the same period of the controllability scheduling.

\item $\Delta$ is the time slot duration.
\end{itemize}
\end{definition}
We assume that each link can be scheduled only one time for each scheduling period. This does not lead to loss of generality, since it is always possible to obtain an equivalent model that satisfies this constraint by appropriately splitting the nodes of the graph, as already illustrated in the memory slot graph definition of \cite{AlurRTAS09}. We define a connectivity property of the controllability and observability graphs with respect to the corresponding scheduling.
\begin{definition}
Given a controllability graph $\mathcal G_{\mathcal R}$ and scheduling $\eta_{\mathcal R}$, we define $\mathcal G_{\mathcal R}(\eta_{\mathcal R}(k))$ the sub-graph of $\mathcal G_{\mathcal R}$ induced by keeping the edges scheduled at time $k$. We define $\mathcal G_{\mathcal R}(\eta_{\mathcal R}) = \bigcup\limits_{k = 1}^{\Pi} \mathcal G_{\mathcal R}(\eta_{\mathcal R}(k))$ the sub-graph of $\mathcal G_{\mathcal R}$ induced by keeping the union of edges scheduled during the whole scheduling period $\Pi$.
\end{definition}

\begin{definition}
We say that a controllability graph $\mathcal G_{\mathcal R}$ is jointly connected by a controllability scheduling $\eta_{\mathcal R}$ if and only if there exists a path from the controller node $v_c$ to the actuator node $v_u$ in $\mathcal G_{\mathcal R}(\eta_{\mathcal R})$. We denote by $D_{\mathcal R}$ the length of the longest path connecting $v_c$ to $v_u$ in $\mathcal G_{\mathcal R}(\eta_{\mathcal R})$.
\end{definition}

The above definitions can be given similarly for observability graph $\mathcal G_{\mathcal O}$ and scheduling $\eta_{\mathcal O}$.

\begin{figure*}[ht]
\begin{center}
\includegraphics[width=0.6\textwidth]{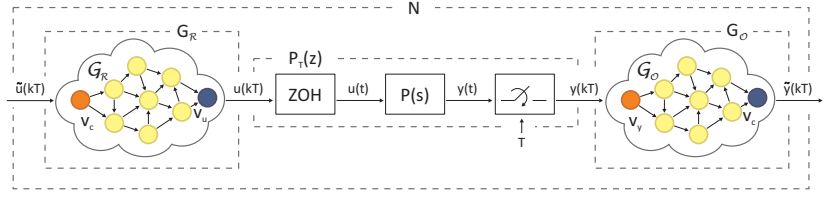}
\caption{MCN interconnected system.} \label{MCNfunctionblocks}
\end{center}
\end{figure*}

The semantics of a MCN $\mathcal N$ can be modeled by the interconnection $N$ of blocks as in Figure \ref{MCNfunctionblocks}. The block $P_{T}$ is characterized by the transfer function $P_{T}(z)$ of the discrete-time plant $P_{T}$, obtained by discretizing $P(s)~=~C(sI~-~A)^{-1}B$ with sampling time ${T} = \Pi \Delta$ equal to the scheduling period duration. The block $G_{\mathcal R}$ models the dynamics introduced by the data flow of the actuation data through the network $\mathcal G_{\mathcal R}$ according to the applied controllability scheduling $\eta_{\mathcal R}$. In order to define the dynamical behavior of $G_{\mathcal R}$, we need to define the semantics of the data flow through the network, according to the scheduling $\eta_{\mathcal R}$.
We associate to the controller node $v_c$ a real value $\mu_c(kT)$ at time $k$, and we assume that $v_c$ is periodically updated with a new control command at the beginning of each scheduling period and holds this value for the whole time duration of the scheduling period. Formally, $\mu_c(kT) = \tilde u(kT)$.

The dynamics of the other nodes needs to be defined at the level of the time slots. We associate to each other node $v_j \in V_{\mathcal R} \setminus \{v_c\}$ a real value $\bar \mu_{i,j}(h)$ at time slot $h$, for each node $v_i$ belonging to the set $inc(v_j) = \{v \in V_{\mathcal R} : (v, v_j) \in E_{\mathcal R}\}$ of edges incoming in $v_j$. Formally:
$$
\mu_{j}(h) = \sum\limits_{v_i \in inc(v_j)} \bar \mu_{i,j}(h)
$$
is the sum of the variables associated to node $v_j$ in the time slot $h$. When the link from $v_i$ to $v_j$ is scheduled at time slot $h$, the variable of node $v_j$ is updated with the sum of the variables of node $v_i$ multiplied for the link weight $W_{\mathcal R}(v_i,v_j)$. Formally, for each $v_j \in V_{\mathcal R} \setminus \{v_c\}$ and for each time slot $h \in \{1, \ldots, \Pi\}$:
\begin{equation*}
\bar \mu_{i,j}(h+1) = \left\{
\begin{array}{ll}
  \bar \mu_{i,j}(h) & \text{if } (v_i, v_j) \notin \eta_{\mathcal R}(h), \\
  W_{\mathcal R}(v_i,v_j) \mu_{i}(h) & \text{if } (v_i, v_j) \in \eta_{\mathcal R}(h).
\end{array}
\right.
\end{equation*}
Finally, the actuator node $v_u$ periodically actuates a new actuation command at the beginning of each scheduling period on the basis of its variable $\mu_u$, and holds this value for the whole time duration of the scheduling period. Formally:
$$
u(kT) = \mu_u(kT) = \sum_{v_i \in inc(v_u)} \bar \mu_{i,u}(kT).
$$
On the basis of the semantics defined above, it is possible to model the dynamical behavior of $G_{\mathcal R}$ as follows.
\begin{proposition}\label{propGz}
Given $\mathcal G_{\mathcal R}$ and $\eta_{\mathcal R}$, the block $G_{\mathcal R}$ can be modeled as a discrete time SISO LTI system with sampling time equal to the scheduling period duration $T = \Pi \Delta$, and characterized by the following transfer function:
\begin{equation*}
G_{\mathcal R}(z) = \sum\limits_{i = 1}^{D_{\mathcal R}} \frac{\gamma_{\mathcal R}(i)}{z^i}, \qquad \gamma_{\mathcal R}(i) = \sum\limits_{\rho \in \chi_{\mathcal R}(i)} W_{\mathcal R}(\rho),
\end{equation*}
where $\forall i \in \{1, \ldots, D_{\mathcal R-1}\}$, $\gamma_{\mathcal R}(i) \in \mathbb R^+ \cup \{0\}$, $\gamma_{\mathcal R}(D_{\mathcal R}) \neq 0$
\end{proposition}
\emph{Proof:} We need to characterize the dynamics of $u(kT) = \mu_u(kT)$ with respect to $\tilde u(kT) = \mu_c(kT)$. Let $\chi_{\mathcal R}$ be the set of all simple paths of $\mathcal G_{\mathcal R}(\eta_{\mathcal R})$ starting from $v_c$ and terminating in $v_u$. We remark that, since $\mathcal G_{\mathcal R}$ is acyclic, then $\chi_{\mathcal R}$ is a finite set. Given any path $\rho = v_c, v_1, \ldots, v_n, v_u  \in \chi_{\mathcal R}$, with $n \in \{0, \ldots, D_{\mathcal R}-1\}$, we define $W_{\mathcal R}(\rho) = W_{\mathcal R}{(v_c, v_1)} W_{\mathcal R}{(v_1, v_2)} \cdots W_{\mathcal R}{(v_{n-1}, v_n)} W_{\mathcal R}{(v_n, v_u)}$ as the product of the weights of the edges of $\rho$. During each scheduling period, each path $\rho \in \chi_{\mathcal R}$ provides an update of the variable of node $v_u$ given by $\mu_u(kT) = W_{\mathcal R}(\rho) \mu_c\left((k-\delta_{\mathcal R}(\rho))T\right)$, where $\delta_{\mathcal R}(\rho)$ is the delay introduced by $\rho$ in terms of scheduling periods. We show how to compute $\delta_{\mathcal R}(\rho)$.

Consider a path $\rho = v_c, v_1, \ldots, v_n, v_u  \in \chi_{\mathcal R}$ with length $|\rho| = n+1$, and assume that the links are scheduled in the same order of the path (as in $\eta'_{\mathcal R}$ of Example \ref{exSinglePath}), i.e. $(v_c, v_1), (v_1, v_2), \ldots, (v_{n-1}, v_n), (v_n, v_u)$. In this case, the information stored in $v_c$ is conveyed to $v_u$ with a delay of just one scheduling period.

Assume now that links are scheduled in the opposite order (as in $\eta''_{\mathcal R}$ of Example \ref{exSinglePath}), i.e. $(v_n, v_u), (v_{n-1}, v_n), \ldots, (v_1, v_2),\allowbreak (v_c, v_1)$. In this case, the information stored in $v_c$ is conveyed to $v_u$ with a delay of $|\rho|$ scheduling periods.

It is easy to show that, given any other scheduling and using the same reasoning, it is possible to determine the delay of scheduling periods $\delta_{\mathcal R}(\rho)$ introduced by the path $\rho$ applying the scheduling $\eta_{\mathcal R}$. It is also easy to verify that $1 \leq \delta_{\mathcal R}(\rho) \leq |\rho| \leq D_{\mathcal R}$. Define $\chi_{\mathcal R}(i) = \{\rho \in \chi_{\mathcal R} : \delta_{\mathcal R}(\rho) = i\}$ the set of all paths that, by applying the scheduling $\eta_{\mathcal R}$, introduce a delay of $i$ scheduling periods to convey the information stored in $v_c$ to $v_u$. We remark that the set $\{\chi_{\mathcal R}(i)\}_{i=1}^{D_{\mathcal R}}$ is a partition of $\chi_{\mathcal R}$.

It is possible to show by induction, starting from the nodes adjacent to $v_c$ up to the node $v_u$, that at the end of a scheduling period each node $v_u$ contains the sum of the contributions of all paths starting from $v_c$ and terminating in $v_u$. Thus, at the end of each scheduling period, the variable associated to the actuator node $v_u$ contains the sum of the contributions of all paths $\rho \in \chi_{\mathcal R}$:
\begin{equation*}
\mu_u(kT) = \sum_{i=1}^{D_{\mathcal R}}\sum_{\rho\in\chi_{\mathcal R}(i)} W_{\mathcal R}(\rho) \mu_c\left((k-i)T\right).
\end{equation*}
Since $u(kT) = \mu_u(kT)$ and $\tilde u(kT) = \mu_c(kT)$, the following holds:
\begin{equation*}
G_{\mathcal R}(z) = \sum\limits_{i = 1}^{D_{\mathcal R}} \frac{\gamma_{\mathcal R}(i)}{z^i}, \qquad \gamma_{\mathcal R}(i) = \sum\limits_{\rho \in \chi_{\mathcal R}(i)} W_{\mathcal R}(\rho).
\end{equation*}
If $\chi_{\mathcal R}(i) = \varnothing$ then $\gamma_{\mathcal R}(i) = 0$, otherwise $\gamma_{\mathcal R}(i) > 0$ because $W_{\mathcal R}$ is positive definite. Thus $\forall i \in \{1, \ldots, D_{\mathcal R}\}, \gamma_{\mathcal R}(i) \in \mathbb R^+ \cup \{0\}$. This completes the proof. $\blacksquare$

\begin{figure}[ht]
\begin{center}
\includegraphics[width=0.5\textwidth]{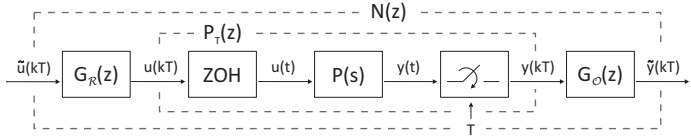}
\caption{Transfer function of the MCN interconnected system.} \label{MCNblocks}
\end{center}
\end{figure}

On the basis of the above reasoning, it is possible to model the semantics of a MCN $\mathcal N$ as in Figure \ref{MCNblocks}, where each block is a discrete time SISO LTI system with sampling time equal to the scheduling period duration, characterized by the transfer functions $G_{\mathcal R}(z)$, $P_T(z)$ and $G_{\mathcal O}(z)$.

The following example motivates the use of redundancy in MCNs characterized by failures in links.
\begin{example}\label{exSinglePath}
Consider a MCN $\mathcal N_1 = ( \mathcal P, \mathcal G_{\mathcal R}, \eta_{\mathcal R}, \mathcal G_{\mathcal O}, \eta_{\mathcal O}, \Delta)$. We remark that, since we are addressing controllability analysis, we do not consider the effect of observability radio connectivity graph and scheduling. $\mathcal P = (A,B,C)$ represents a continuous-time SISO LTI plant characterized by the following transfer function:
\begin{equation*}
P(s) = \frac{1}{(s^2 - 69.31s + 25880) (s^2 + 24670)}.
\end{equation*}
We consider a single-path scenario for the controllability radio connectivity graph $\mathcal G_{\mathcal R}$ as shown in Figure \ref{singlePath}, where $v_c = v_1$ and $v_u = v_3$.
\begin{figure}[b]
\begin{center}
\includegraphics[width=0.5\textwidth]{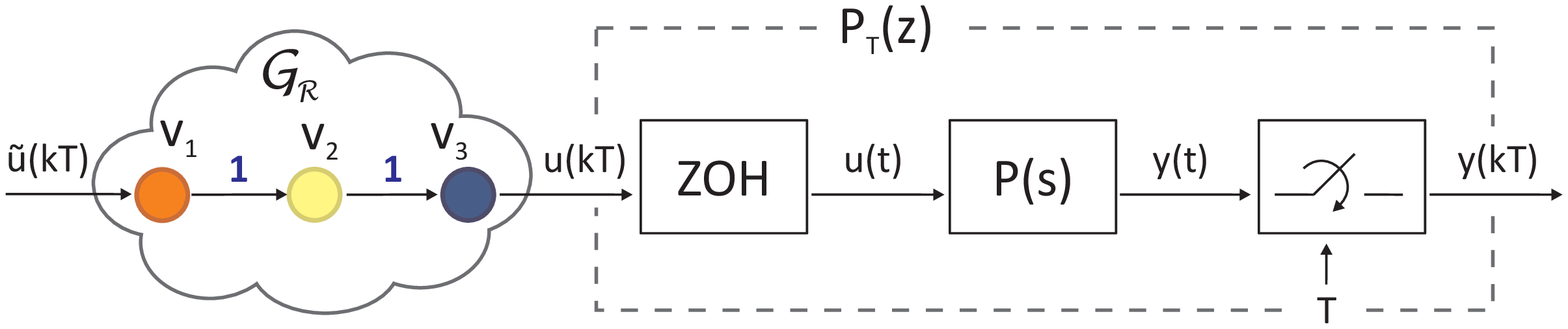}
\caption{Single-path scenario graph $\mathcal G_{\mathcal R}$.} \label{singlePath}
\end{center}
\end{figure}
Let us define two schedules on $\mathcal G_{\mathcal R}$, which convey actuation data from $v_c$ to $v_u$: we define $\eta'_{\mathcal R}$ as the string $\langle\{(v_1, v_2)\},\{(v_2, v_3)\}\rangle$, and $\eta''_{\mathcal R}$ as the string $\langle\{(v_2, v_3)\},\{(v_1, v_2)\}\rangle$. Both scheduling have period $\Pi = 2$. Let $\Delta$ be the duration in seconds of a single time slot of the scheduling, for example $10\ ms$ as in WirelessHART: then the duration of the whole scheduling period is given by $T = \Pi \Delta = 2 \cdot 10\ ms = 20\ ms$.
\begin{figure}[t]
\begin{center}
\includegraphics[width=0.3\textwidth]{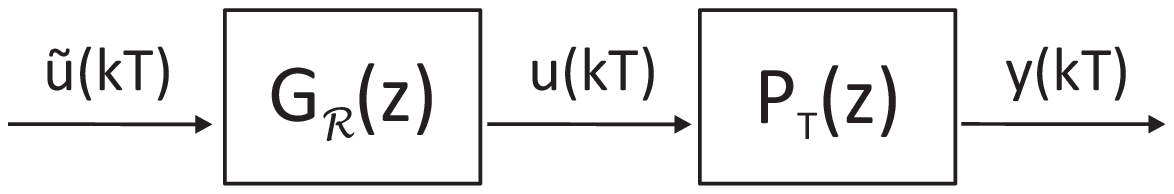}
\caption{Single-path scenario interconnected system.} \label{singlePathScheme}
\end{center}
\end{figure}
If we apply $\eta'_{\mathcal R}$, the block $G_{\mathcal R}$ introduces a delay equal to $1$ scheduling period, i.e. $u(kT)= \tilde u((k-1)T)$. In fact, using $\eta'_{\mathcal R}$, the value of $\tilde u$ is conveyed first to $v_2$ and then to $v_3$ in just one scheduling period. Thus, the block $G_{\mathcal R}$ is characterized by the transfer function $G_{\mathcal R}(z) = \frac{1}{z}$.

If we apply $\eta''_{\mathcal R}$, the block $G_{\mathcal R}$ introduces a delay equal to~$2$ scheduling periods, i.e. $u(kT)=\tilde u((k-2)T)$. In fact, using $\eta''_{\mathcal R}$, the value of $\tilde u$ is conveyed to $v_2$ in the first period, and to $v_3$ in the second period. Thus, the block $G_{\mathcal R}$ is characterized by the transfer function $G_{\mathcal R}(z) = \frac{1}{z^2}$.

It is easy to verify that, in both cases, the cascade~$M$ of systems $G_{\mathcal R}$ and $P_T$ shown in Figure~\ref{singlePathScheme} always satisfies the controllability condition because $G_{\mathcal R}(z)$ does not have any zeros, so its interconnection to $P_T(z)$ can not introduce any pole cancelation. However, if one of the links of $\mathcal G_{\mathcal R}$ is damaged, then $G_{\mathcal R}(z) = 0$ and the MCN does not satisfy the controllability condition. Moreover, an error in one of the data transmissions between $v_1$, $v_2$ and $v_3$ is totally transferred to the input of the plant. $\square$
\end{example}
The above example motivates the exploitation of redundancy, for instance by sending control data through multiple paths in the same scheduling period and then merging these components in the actuator node. We call this approach redundancy by static multi-path routing. An alternative is sending control data through a single-path (route) for each scheduling period, and dynamically updating this route in order to avoid faulty nodes in the new route. We call this approach redundancy by dynamic single-path routing. At the best of our knowledge, although there exist several algorithms for static (e.g. Dijkstra and Bellman-Ford) and dynamic routing \cite{DynamicRouting} of multi-hop networks, none of them have been designed to address control specifications. In this paper, we only address the problem of designing redundancy by static multi-path routing, in order to preserve controllability and observability structural properties of a MCN.

The \emph{Pros} of applying redundancy by static multi-path routing are the following: first, controllability of the MCN is robust to failures of links; and second, the effect of data transmission errors on a single link is alleviated when averaging all the components received from the multiple paths. It is worth to remark that, although protocols such as WirelessHART and ISA-100 are oriented to single-path routing (i.e. sensing and actuation data are sent to the controller via a unique path of wireless nodes), it is possible to implement redundancy of sensing and actuation data by appropriately defining the scheduling in order to achieve multi-path routing.

The \emph{Cons} of redundancy by static multi-path routing are the following: first, it increases data traffic in the network, but this is a necessary investment to improve robustness with respect to link failures when we apply static routing; second, sending control data through multiple paths and then merging them in the actuator node generates dynamics that might invalidate controllability conditions, as illustrated in the following example.
\begin{example}\label{exRedundancy}
Let define a MCN $\mathcal N_2 = ( \mathcal P, \mathcal G_{\mathcal R}, \eta_{\mathcal R}, \mathcal G_{\mathcal O}, \eta_{\mathcal O}, \Delta)$. As above, we do not consider the observability radio connectivity graph and scheduling. $\mathcal P = (A,B,C)$ represents a continuous-time SISO LTI plant characterized by the transfer function $P(s)$ adopted in the previous example. We consider a multi-path scenario for the controllability radio connectivity graph $\mathcal G_{\mathcal R}$ as shown in Figure \ref{multiPath} where $v_c = v_1$ and $v_u = v_7$.
\begin{figure}[t]
\begin{center}
\includegraphics[width=0.5\textwidth]{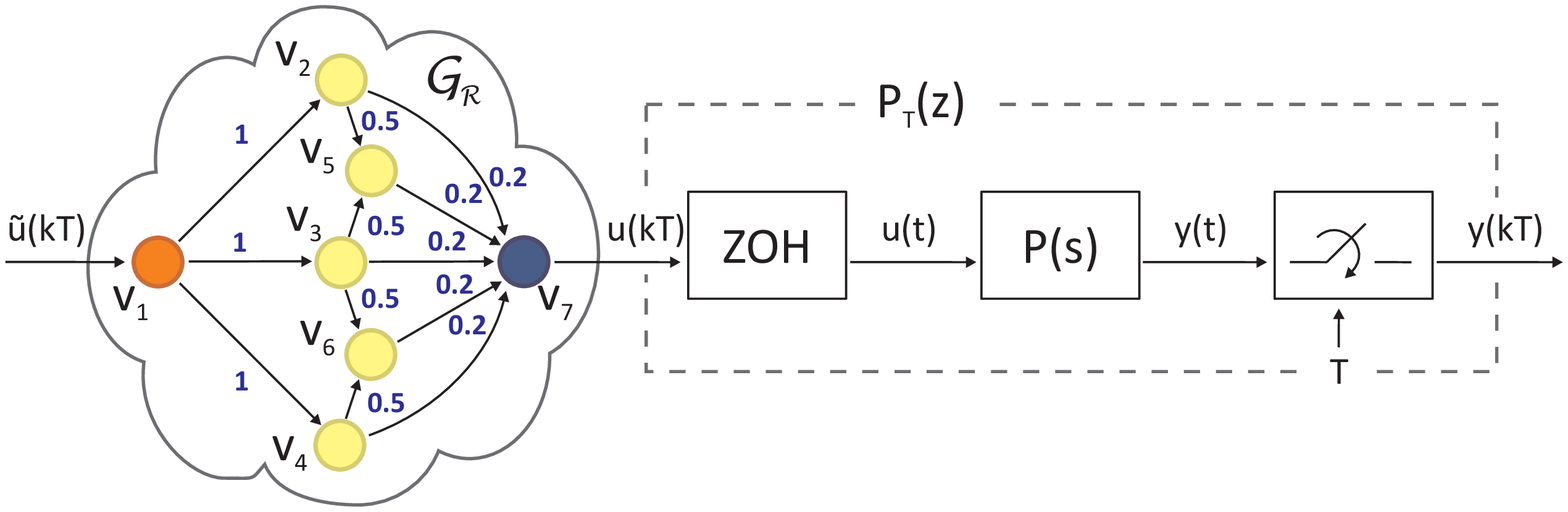}
\caption{Multi-path scenario graph $\mathcal G_{\mathcal R}$.} \label{multiPath}
\end{center}
\end{figure}
Differently from the controllability graph of Example \ref{exSinglePath}, in this case some nodes receive actuation data from multiple links. We define the weight function $W_{\mathcal R}$ so that all nodes equally weight the contribution of each incoming link.

Let us define three schedules on $\mathcal G_{\mathcal R}$, which convey actuation data from $v_c$ to $v_u$ using multiple paths:
\begin{align*}
\eta_{\mathcal R}^{a} = \langle& \{(v_1, v_2), (v_1, v_3), (v_1, v_4), (v_2, v_5), (v_2, v_7), (v_3, v_5),\\
&(v_3, v_6), (v_3, v_7), (v_4, v_6), (v_4, v_7), (v_5, v_7), (v_6, v_7)\}\rangle,\\
\eta_{\mathcal R}^{b} = \langle& \{(v_1, v_2), (v_1, v_3), (v_1, v_4)\}, \{(v_2, v_5), (v_2, v_7)\},\\
&\{(v_3, v_5), (v_3, v_6), (v_3, v_7)\},\{(v_4, v_6), (v_4, v_7)\},\\
&\{(v_5,v_7)\},\{(v_6,v_7)\}\rangle,\\
\eta_{\mathcal R}^{c} = \langle&\{(v_1, v_2), (v_1, v_3), (v_1, v_4), (v_5, v_7), (v_6, v_7)\}, \{(v_2, v_5),\\
&(v_2, v_7), (v_3, v_5), (v_3, v_6), (v_3,v_7), (v_4, v_6), (v_4, v_7)\}\rangle.
\end{align*}
The scheduling have periods $\Pi^a = 1$, $\Pi^b = 6$, and $\Pi^c = 2$. Let $\Delta = 10\ ms$ as above: then the durations of the scheduling periods are given by $T^a = 10\ ms$, $T^b = 60\ ms$, and $T^c = 20\ ms$.

Scheduling $a$ consists of only one time slot, where all nodes transmit simultaneously. This is a corner case, that is realistic e.g. if the MAC layer of the communication protocol implements CDMA. The scheduling $a$ produces the following dynamics for $u(kT^a)$:
\begin{equation*}
u(kT^a) = \frac{3}{5}\tilde{u}((k-2)T^a)+\frac{2}{5}\tilde{u}((k-3)T^a).
\end{equation*}
Note that a delay of $2$ scheduling periods of $\eta_{\mathcal R}^{a}$ corresponds to an actuation delay of $2 \cdot T^a = 20\ ms$.

Scheduling $b$ consists of $6$ time slots: in each time slot only one node transmits and the other nodes receive. This is also a corner case, that is realistic when the communication protocol implements TDMA (i.e. only one node is allowed to transmit for each time slot). The scheduling $b$ produces the following dynamics for $u(kT^b)$:
\begin{equation*}
u(kT^b) = \tilde{u}((k-1)T^b).
\end{equation*}
Note that a delay of $1$ scheduling period of $\eta_{\mathcal R}^{b}$ corresponds to an actuation delay of $1 \cdot T^b = 60\ ms$.

Scheduling $c$ consists of $2$ time slots: in the first time slot only nodes $v_1$, $v_5$ and $v_6$ are allowed to transmit, since they are assumed not to interfere each other; in the second time slot only nodes $v_2$, $v_3$ and $v_4$ are allowed to transmit. This scheduling is a tradeoff between scheduling $a$ and $b$, and is realistic when the communication protocol implements mixed TDMA and FDMA. The scheduling $c$ produces the following dynamics for $u(kT^c)$:
\begin{equation*}
u(kT^c) = \frac{3}{5}\tilde{u}((k-1)T^c)+\frac{2}{5}\tilde{u}((k-2)T^c).
\end{equation*}
Note that a delay of 1 scheduling period $\eta_{\mathcal R}^{c}$ corresponds to an actuation delay of $1 \cdot T^c = 20\ ms$.

As a comparison to the single-path scenario of Example~\ref{exSinglePath}, we can conclude that adding redundancy generally increases the actuation delay of the control input, and thus worsen responsiveness of the control algorithm. The only case when the actuation delay does not increase is when we dispose of a multi-hop network that allows simultaneous transmission of all links (scheduling $a$, e.g. using CDMA, $T^a = 10\ ms$). Unfortunately, this is not generally the case for current specifications for wireless networks, e.g. WirelessHART and ISA-100 do not admit CDMA.

In the worst scenario, only one node can transmit data for each time slot (scheduling $b$, e.g. using TDMA, $T^b~=~60\ ms$). In this case, the actuation delay strongly increases, in contrast with the requirement of designing responsive control algorithms. In typical industrial scenarios, as illustrated in \cite{D'InnocenzoCASE2009}, if we do not admit simultaneous transmission of two nodes then it is not always possible to design a scheduling that satisfies constraints on the actuation delay guaranteeing the achievement of control specifications on the closed loop system.

For this reason, we choose a scheduling that allows multiple transmissions on a subset of links, which do not interfere each other (scheduling $c$, e.g. using TDMA and FDMA, $T^c = 20\ ms$). The schedule $c$ has the advantage of introducing redundancy and only moderately increasing the actuation delay. Moreover, it is reasonably implementable since it requires simultaneous use of mixed TDMA and FDMA, which is standardized in existing communication protocols for wireless sensor networks such as WirelessHART and ISA-100.

Using $\eta_{\mathcal R}^{c}$, the system can be seen as the cascade $M^c$ of blocks as in Figure \ref{multiPathScheme}. The block $P_{T^c}$ is characterized by the transfer function $P_{T^c}(z)$ obtained by discretizing $P(s)$ with sampling time ${T^c}$ as follows:
\begin{equation*}
P_{T^c}(z) = 4.2932\times 10^{-9} \frac{z + 1.189}{(z + 1)(z + 2)},
\end{equation*}
where the poles of $P_{T^c}(z)$ are $p_1=-1$ and $p_2=-2$. The block $G_{\mathcal R}$ models the dynamics introduced by the data flow through the network $\mathcal G_{\mathcal R}$, according to the applied scheduling $\eta^c_{\mathcal R}$.
\begin{figure}[ht]
\begin{center}
\includegraphics[width=0.3\textwidth]{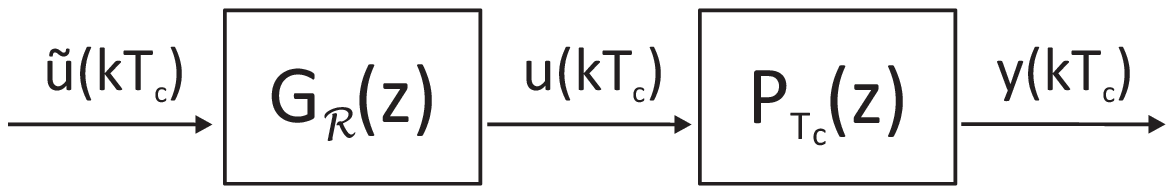}
\caption{Multi-path scenario interconnected system.}\label{multiPathScheme}
\end{center}
\end{figure}

The block $G_{\mathcal R}$ is characterized by the transfer function
\begin{equation*}
G_{\mathcal R}(z) = \frac{\frac{3}{5} z + \frac{2}{5}}{z^2} \Rightarrow z=-\frac{2}{3}.
\end{equation*}
Thus, when no failures occur, the zero $z$ of $G_{\mathcal R}(z)$ does not coincide with any of the poles of $P_{T^c}(z)$. Assume that a failure occurs in the link $(v_3,v_7)$, denoted by the failures configuration $f_1$. Then the following holds:
\begin{equation*}
f_1 = \left\{(v_3,v_7)\right\} \Rightarrow G_{\mathcal R}^{f_1}(z) = \frac{\frac{2}{5}z+\frac{2}{5}}{z^2} \Rightarrow z_{f_1}=-1.
\end{equation*}
Thus, when the fault $f_1$ occurs the zero $z_{f_1}$ of $G^{f_1}_{\mathcal R}(z)$ coincides with the pole $p_1$ of $P_{T^c}(z)$. Assume that a failure simultaneously occurs in the links $(v_3,v_7)$ and $(v_5,v_7)$, denoted by the failures configuration $f_2$. Then the following holds:
\begin{equation*}
f_2 = \left\{(v_3,v_7),(v_5,v_7)\right\} \Rightarrow G_{\mathcal R}^{f_2}(z) = \frac{\frac{2}{5}z+\frac{1}{5}}{z^2} \Rightarrow z_{f_2}=-\frac{1}{2}.
\end{equation*}
Thus, when the fault $f_2$ occurs, the zero $z_{f_2}$ of $G^{f_2}_{\mathcal R}(z)$ does not coincide with any of the poles of $P_{T^c}(z)$. Assume that a failure simultaneously occurs in the links $(v_3,v_7)$ and $(v_2,v_7)$, denoted by the failures configuration $f_3$. Then the following holds:
\begin{equation*}
f_3 = \left\{(v_3,v_7),(v_2,v_7)\right\} \Rightarrow G_{\mathcal R}^{f_3}(z) = \frac{\frac{1}{5}z+\frac{2}{5}}{z^2} \Rightarrow z_{f_3}=-2.
\end{equation*}
Thus, when the fault $f_3$ occurs the zero $z_{f_3}$ of $G^{f_3}_{\mathcal R}(z)$ coincides with the pole $p_2$ of $P_{T^c}(z)$. Note that, according to a generic set of faulty edges $f$ and to the weight function $W_{\mathcal R}$, some failures introduce dynamics in the block $G^{f}_{\mathcal R}(z)$ that invalidate controllability of $P_{T^c}(z)$. $\square$
\end{example}
Example \ref{exRedundancy} motivates characterizing MCNs controllability conditions over the plant dynamics, the scheduling and the routing. We will address this problem in the next Section.


\section{Stabilizability of MCNs} \label{secMCNStabilizability}

As discussed in Section \ref{secMCNModeling}, MCNs $\mathcal N$ can be modeled by the interconnection $N$ of blocks as in Figure \ref{MCNblocks}.
\begin{definition}
We say that a MCN $\mathcal N$ is controllable (resp. observable) if and only if $N$ is controllable (resp. observable). Moreover, we say that $\mathcal N$ is stabilizable (resp. detectable) if and only if $N$ is stabilizable (resp. detectable).
\end{definition}
\begin{theorem}\label{thMCNControllability}
A MCN $\mathcal N$ is controllable if and only if the following hold:
\begin{enumerate}
\item $(A,B)$ is controllable;

\item $\mathcal G_{\mathcal R}$ is jointly connected by $\eta_{\mathcal R}$;

\item for each pole $p$ of $P_T(z)$, $\sum\limits_{i=1}^{D_{\mathcal R}} \gamma_{\mathcal R}(i) p^{D_{\mathcal R}-i} \neq 0$;

\item for each zero $z$ of $P_T(z)$, $z \neq 0$.
\end{enumerate}
\emph{Proof:} (\textbf{Sufficiency}) Since Condition 1 states that the block $P_T$ is controllable, we need to prove that Conditions 2,3 and 4 imply that $G_{\mathcal R}(z) \neq 0$, and that poles of $P_T(z)$ and $G_{\mathcal O}(z)$ do not coincide with zeros of $G_{\mathcal R}(z)$ and $P_T(z)$.

We prove that Condition 2 implies that $G_{\mathcal R}(z) \neq 0$. If $\mathcal G_{\mathcal R}$ is jointly connected by $\eta_{\mathcal R}$, then there exists at least one path $\rho = v_c, \ldots, v_u$ with $W_{\mathcal R}(\rho) > 0$. Let $\delta_{\mathcal R}(\rho) = i$, with $1 \leq i \leq |\rho| \leq D_{\mathcal R}$. Since $W_{\mathcal R}$ is positive definite, then $\gamma_{\mathcal R}(i) > W_{\mathcal R}(\rho) > 0$. This implies that $G_{\mathcal R}(z) \neq 0$.

We prove that Condition 3 implies that any zero of $G_{\mathcal R}(z)$ does not coincide with any pole of $P_T(z)$. We can write $G_{\mathcal R}(z)$ as follows:
\begin{equation*}
G_{\mathcal R}(z) = \frac{\sum\limits_{i=1}^{D_{\mathcal R}} \gamma_{\mathcal R}(i) z^{D_{\mathcal R} -i}}{z^D_{\mathcal R}}.
\end{equation*}
Thus, $p$ is a zero of $G_{\mathcal R}(z)$ if and only if Condition 3 holds.

Since $G_{\mathcal R}(z)$ can not contain zeros in 0 and $G_{\mathcal O}(z)$ only has poles in 0, then any zero of $G_{\mathcal R}(z)$ does not coincide with any pole of $G_{\mathcal O}(z)$.

Since $G_{\mathcal R}(z)$ only has poles in 0, Condition 4 implies that any zero of $P_T(z)$ does not coincide with any pole of $G_{\mathcal R}(z)$.

This completes the first part of the proof.

\smallskip

(\textbf{Necessity}) For each condition, we assume that it is not satisfied, and prove that this implies that $N$ is not controllable.

Assume that Condition 1 is not satisfied, then clearly the system $N$ is not controllable.

Assume that Condition 2 is not satisfied. Since $\chi_{\mathcal R} = \varnothing$, then $\gamma_{\mathcal R}(i) = 0$ for each $i \in \{1, \ldots D_{\mathcal R}\}$, and thus $G_{\mathcal R}(z)=0$. This clearly implies that the system $N$ is not controllable.

Assume that Condition 3 is not satisfied, then there exists a zero of $G_{\mathcal R}(z)$ that coincides with a pole of $P_T(z)$. This clearly implies that the system $N$ is not controllable.

Assume that Condition 4 is not satisfied, then there exists a zero of $P_T(z)$ that coincides with a pole of $G_{\mathcal O}(z)$. This clearly implies that the system $N$ is not controllable. $\blacksquare$
\end{theorem}
Note that controllability of $(A,B)$ and connectivity of the controller and actuator nodes are necessary conditions, as suggested by the intuition, but they are not sufficient. In fact, the controllability scheduling may introduce dynamics that invalidate controllability. This issue generates Conditions 3 and 4 of Theorem \ref{thMCNControllability}, that provide together with Conditions 1 and 2 necessary and sufficient controllability conditions.

Another interesting remark is that, in order to guarantee controllability, we do not need to design an ordered schedule, i.e. we can schedule links with any order. This is an interesting result, since it allows much more freedom in designing the scheduling. The main problem, as illustrated in Example \ref{exRedundancy}, is the design of a weight function $W_{\mathcal R}$ such that Condition 3 of Theorem \ref{thMCNControllability} is satisfied: we address and solve this issue in the following Section. Note that designing both a scheduling function $\eta_{\mathcal R}$ and a weight function $W_{\mathcal R}$ that satisfy the conditions of the Theorem~\ref{thMCNControllability}, corresponds to designing a scheduling and a multi-path routing of the communication protocol.

The following corollary can be proved using the same reasoning as in the proof of Theorem \ref{thMCNControllability}.
\begin{corollary}\label{thMCNStabilizability}
A MCN $\mathcal N$ is stabilizable if and only if the following hold:
\begin{enumerate}
\item $(A,B)$ is stabilizable;

\item $\mathcal G_{\mathcal R}$ is jointly connected by $\eta_{\mathcal R}$;

\item for each pole $p$ of $P_T(z)$ such that $|p| \geq 1$, $\sum\limits_{i=1}^{D_{\mathcal R}} \gamma_{\mathcal R}(i) p^{D_{\mathcal R}-i} \neq 0$.
\end{enumerate}
\end{corollary}
By duality:
\begin{corollary}\label{thMCNObservability}
A MCN $\mathcal N$ is observable if and only if the following hold:
\begin{enumerate}
\item $(C,A)$ is observable;

\item $\mathcal G_{\mathcal O}$ is jointly connected by $\eta_{\mathcal O}$;

\item for each pole $p$ of $P_T(z)$, $\sum\limits_{i=1}^{D_{\mathcal O}} \gamma_{\mathcal O}(i) p^{D_{\mathcal O}-i} \neq 0$;

\item for each zero $z$ of $P_T(z)$, $z \neq 0$.
\end{enumerate}
\end{corollary}
\begin{corollary}\label{thMCNDetectablity}
A MCN $\mathcal N$ is detectable if and only if the following hold:
\begin{enumerate}
\item $(C,A)$ is detectable;

\item $\mathcal G_{\mathcal O}$ is jointly connected by $\eta_{\mathcal O}$;

\item for each pole $p$ of $P_T(z)$ such that $|p| \geq 1$, $\sum\limits_{i=1}^{D_{\mathcal O}} \gamma_{\mathcal O}(i) p^{D_{\mathcal O}-i} \neq 0$.
\end{enumerate}
\end{corollary}


\section{Fault tolerant Stabilizability of MCNs} \label{secMCNRobustStabilizability}

Given a MCN and a set $f \subseteq E_{\mathcal R} \cup E_{\mathcal O}$ of communication links subject to a failure, we define a faulty MCN as follows:
\begin{definition}
Given MCN $\mathcal N = ( \mathcal P, \mathcal G_{\mathcal R}, \eta_{\mathcal R}, \mathcal G_{\mathcal O}, \eta_{\mathcal O}, \Delta)$, let $f \subseteq E_{\mathcal R} \cup E_{\mathcal O}$ be a set of faulty links. We define the faulty MCN $\mathcal N_f = ( \mathcal P, \mathcal G_{\mathcal R}, \eta^f_{\mathcal R}, \mathcal G_{\mathcal O}, \eta^f_{\mathcal O}, \Delta)$, where $\forall k~\in~\{1, \ldots, \Pi\}$, $\eta^f_{\mathcal R}(k) = \eta_{\mathcal R}(k) \setminus(\eta_{\mathcal R}(k) \cap f)$ and $\eta^f_{\mathcal O}(k) = \eta_{\mathcal O}(k) \setminus (\eta_{\mathcal O}(k) \cap f)$.
\end{definition}
In other words, the faulty MCN $\mathcal N_f$ is obtained by removing the faulty links from the schedules, while keeping the original radio connectivity graphs and the weight functions.

Let $\mathcal F \subseteq 2^{E_{\mathcal R} \cup E_{\mathcal O}}$ be a set of failures configurations. The empty set $\varnothing$ always belong to $\mathcal F$, and $\mathcal N_{\varnothing}$ represents the MCN in absence of failures. As clearly illustrated in Section \ref{secMCNStabilizability}, the main issue in designing a MCN such that it is controllable even with link failures, is the choice of the weight function $W_{\mathcal R}$. Since we are exploiting redundancy by static multi-path routing, we need to design a unique static weight function $\bar W_{\mathcal R}$ (which implicitly defines the weight of each routing path) such that $\mathcal N_f$ is controllable for each $f \in \mathcal F$. This problem is not trivial: as an example, notice that it is not always possible to arbitrarily assign the value $W_{\mathcal R}(\rho)$ for each path of an acyclic graph that consists of more than 5 vertices, since the number of paths, which corresponds to the number of constraints, is greater than the number of edges, which corresponds to the number of free variables.

The following Theorem provides necessary and sufficient conditions for guaranteeing, given a MCN and a set of failures configurations $\mathcal F$, the existence of a weight function $\bar W_{\mathcal R}$ such that $\mathcal N_f$ is controllable for each $f \in \mathcal F$. The proof is constructive, and thus provides an algorithm to design $\bar W_{\mathcal R}$.

\begin{theorem}\label{thfaultyMCNControllability}
Given a MCN $\mathcal N$ and a faulty set $\mathcal F$, there exists a weight function $\bar W_{\mathcal R}$ such that $\mathcal N_f$ is controllable for each $f \in \mathcal F$, if and only if the following hold:
\begin{enumerate}
\item $(A,B)$ is controllable;

\item for each $f \in \mathcal F$, $\mathcal G_{\mathcal R}$ is jointly connected by $\eta^f_{\mathcal R}$;

\item for each zero $z$ of $P_T(z)$, $z \neq 0$.
\end{enumerate}
\emph{Proof:} The necessity is trivial, since if one of the above conditions is false then the MCN is not controllable for each weight function $\bar W_{\mathcal R}$. To prove the sufficiency, and since Conditions 1,2 and 4 of Theorem \ref{thMCNControllability} are already satisfied by assumption, we need to provide an algorithm to design $\bar W_{\mathcal R}$.

Pick any weight function $W_{\mathcal R}$: e.g. we can pick a weight function as in Example \ref{exRedundancy} that equally weights all incoming edges to any vertex, i.e. such that for any vertex $v$ of $\mathcal G_{\mathcal R}$ the following holds:
\begin{align*}
\forall v' \in inc(v), W_{\mathcal R}(v',v) = \frac{1}{|inc(v)|}.
\end{align*}
Given $\mathcal G_{\mathcal R} = (V_{\mathcal R}, E_{\mathcal R}, W_{\mathcal R})$, $\eta_{\mathcal R}$ and $\mathcal F$, we define a non-empty set $\bar E_{\mathcal R} = \{\bar e_1, \ldots, \bar e_m\} \subseteq E_{\mathcal R}$ and a partition $\{\mathcal F_i\}_{i=1}^{m}$ of $\mathcal F$ as follows.

Pick $f \in \mathcal F$ and $\bar e_1 \in E_{\mathcal R} \setminus f$. Assign $\bar E_{\mathcal R} = \{\bar e_1\}$ and define the set $\mathcal F_1 = \{f \in \mathcal F : \bar e_1 \notin f\}$. If $\mathcal F_1 = \mathcal F$ we stop the algorithm. Otherwise proceed with the algorithm applying the following inductive $K$-th step: pick $f \in \mathcal F \setminus \bigcup\limits_{i=1}^{K-1} \mathcal F_i$ and $\bar e_K \in E_{\mathcal R} \setminus (\bar E_{\mathcal R} \cup f)$. Assign $\bar E_{\mathcal R} = \bar E_{\mathcal R} \cup \{\bar e_K\}$ and define the set $\mathcal F_K = \{f \in \mathcal F \setminus \bigcup\limits_{i=1}^{K-1} \mathcal F_i : \bar e_K \notin f\}$. If $\bigcup\limits_{i=1}^{K-1} \mathcal F_i \cup \mathcal F_K = \mathcal F$ we stop the algorithm. By iteration, we construct a set $\bar E_{\mathcal R} = \{\bar e_1, \ldots, \bar e_m\} \subset E_{\mathcal R}$ such that the following hold:
\begin{enumerate}
\item the set $\{\mathcal F_i\}_{i=1}^{m}$ is a partition of $\mathcal F$, with $1 \leq m \leq |\mathcal F|$;

\item for each $\bar e_i$ and for each $f \in \bigcup\limits_{j=i}^{m} \mathcal F_j$, then $\bar e_i \in f$.
\end{enumerate}
Note that, since $\mathcal G_{\mathcal R}$ is jointly connected by $\eta^f_{\mathcal R}$, then for any iteration $i$, $\mathcal F_i$ has at least cardinality 1, and thus the algorithm terminates in a finite number of steps upper bounded by $|\mathcal F|$.

Define a weight function $W_{\mathcal R, (\varepsilon_1, \ldots, \varepsilon_m)}$ where $\forall i = \{1, \ldots, m\}, W_{\mathcal R,(\varepsilon_1, \ldots, \varepsilon_m)}(\bar e_i) = W_{\mathcal R}(\bar e_i) + \varepsilon_i$, $\varepsilon_i \in \mathbb R^+$.

Let $\gamma_{\mathcal R, f}(i)$ be given by constructing (as in Proposition \ref{propGz}) the transfer function $\mathcal G^f_{\mathcal R}(z)$ from $\mathcal G_{\mathcal R}$ and $\eta^f_{\mathcal R}$. Condition 3 of Theorem \ref{thMCNControllability} requires that:
\begin{align}\label{eqConditionFaulty_m}
&\forall p \in \mathcal P, \forall f_m \in \mathcal F_m,\notag\\
&\sum\limits_{i = 1}^{D_{\mathcal R}} \gamma_{\mathcal R, f_m}(i) p^{D_{\mathcal R}-i} + b_{p, f_m} \varepsilon_{m} = a_{p, f_m} + b_{p, f_m} \varepsilon_m \neq 0,
\end{align}
where $\mathcal P$ is the set of poles of $P_T(z)$ and $a_{p, f_m}$, $b_{p, f_m}$ are real constants. Note that for all $i \in \{1, \ldots, m-1\}$, $\varepsilon_i$ does not appear in inequalities \eqref{eqConditionFaulty_m} since for each $f_m \in \mathcal F_m$, $\bar e_i \in f_m$ and thus the associated weight $W_{\mathcal R}(\bar e_i)$ does not appear in any path. Pick a value:
\begin{equation*}
\bar \varepsilon_{m} \in \mathbb R^+ \setminus \bigcup\limits_{p \in \mathcal P,\\ f_{m} \in \mathcal F_{m}}\{-\frac{a_{p, f_{m}}}{b_{p,f_{m}}}\},
\end{equation*}
i.e. such that inequalities \eqref{eqConditionFaulty_m} are satisfied, and consider $W_{\mathcal R, (\varepsilon_1, \ldots, \varepsilon_{m-1}, \bar \varepsilon_m)}$. Condition 3 of Theorem \ref{thMCNControllability} requires that:
\begin{align}\label{eqConditionFaulty_m-1}
&\forall p \in \mathcal P, \forall f_{m-1} \in \mathcal F_{m-1},\notag\\
&\sum\limits_{i = 1}^{D_{\mathcal R}} \gamma_{\mathcal R, f_{m-1}}(i) p^{D_{\mathcal R}-i} + b_{p, f_{m-1}} \varepsilon_{m-1} + c_{m,p, f_{m-1}} \bar\varepsilon_{m} \notag\\
=&a_{p, f_{m-1}} + b_{p, f_{m-1}} \varepsilon_{m-1} \neq 0,
\end{align}
where $a_{p, f{m-1}}$, $b_{p, f_{m-1}}$ and $c_{m,p, f_{m-1}}$ are real constants. Note that for all $i \in \{1, \ldots, m-2\}$, $\varepsilon_i$ does not appear in inequalities \eqref{eqConditionFaulty_m-1} since for each $f_{m-1} \in \mathcal F_{m-1}$, $\bar e_i \in f_{m-1}$ and thus the associated weight $W_{\mathcal R}(\bar e_i)$ does not appear in any path. Pick a value:
\begin{equation*}
\bar \varepsilon_{m-1} \in \mathbb R^+ \setminus \bigcup\limits_{p \in \mathcal P,\\ f_{m-1} \in \mathcal F_{m-1}}\{-\frac{a_{p, f_{m-1}}}{b_{p,f_{m-1}}}\},
\end{equation*}
i.e. such that inequalities \eqref{eqConditionFaulty_m-1} are satisfied, and consider $W_{\mathcal R,(\varepsilon_1, \ldots, \varepsilon_{m-2}, \bar \varepsilon_{m-1}, \bar \varepsilon_m)}$. By iteration, after $m-1$ steps we define $W_{\mathcal R,(\varepsilon_1, \bar \varepsilon_2, \ldots, \bar \varepsilon_m)}$. Condition 3 of Theorem \ref{thMCNControllability} requires that:
\begin{align}\label{eqConditionFaulty_1}
&\forall p \in \mathcal P, \forall f_{1} \in \mathcal F_{1},\notag\\
&\sum\limits_{i = 1}^{D_{\mathcal R}} \gamma_{\mathcal R, f_{1}}(i) p^{D_{\mathcal R}-i} + b_{p, f_{1}} \varepsilon_{1} + c_{2,p, f_{1}} \bar\varepsilon_{2} + \ldots + c_{m,p, f_{1}} \bar\varepsilon_{m}\notag\\
=& a_{p, f_{1}} + b_{p, f_{1}} \varepsilon_{1} \neq 0,
\end{align}
where $a_{p, f{1}}$, $b_{p, f_{1}}$ and $c_{2,p, f_{1}}, \ldots, c_{m, p, f_{1}}$ are real constants. Pick a value:
\begin{equation*}
\bar \varepsilon_{1} \in \mathbb R^+ \setminus \cup_{p}\{-\frac{a_{p, f_{1}}}{b_{f_{1}}}\},
\end{equation*}
i.e. such that inequalities \eqref{eqConditionFaulty_1} are satisfied, and consider $W_{\mathcal R,(\bar \varepsilon_1, \ldots, \bar \varepsilon_{m-1}, \bar \varepsilon_m)}$. It is clear that, for each $f \in \bigcup\limits_{i=1}^{m} \mathcal F_i = \mathcal F$, $\mathcal N_f$ satisfies Condition 3 of Theorem \ref{thMCNControllability}. This implies that for each $f \in \mathcal F$, $\mathcal N_f$ is controllable. Assigning $\bar W_{\mathcal R} = W_{\mathcal R,(\bar \varepsilon_1, \ldots, \bar \varepsilon_{m-1}, \bar \varepsilon_m)}$ concludes the proof. $\blacksquare$
\end{theorem}
In the particular case $\mathcal F = \{\varnothing\}$, Theorem \ref{thfaultyMCNControllability} provides necessary and sufficient conditions for guaranteeing, given a non-faulty MCN $\mathcal N_{\varnothing}$, the existence of a weight function $\bar W_{\mathcal R}$ such that $\mathcal N_{\varnothing}$ is controllable, thus solving the design problem defined in the previous Section.

The following corollary can be proved using the same reasoning as in the proof of Theorem \ref{thfaultyMCNControllability}.
\begin{corollary}
Given a MCN $\mathcal N$ and a faulty set $\mathcal F$, there exists a weight function $\bar W_{\mathcal R}$ such that $\mathcal N_f$ is stabilizable for each $f \in \mathcal F$, if and only if the following hold:
\begin{enumerate}
\item $(A,B)$ is stabilizable;

\item for each $f \in \mathcal F$, $\mathcal G_{\mathcal R}$ is jointly connected by $\eta^f_{\mathcal R}$.
\end{enumerate}
\end{corollary}
By duality:
\begin{corollary}
Given a MCN $\mathcal N$ and a faulty set $\mathcal F$, there exists a weight function $\bar W_{\mathcal O}$ such that $\mathcal N_f$ is observable for each $f \in \mathcal F$, if and only if the following hold:
\begin{enumerate}
\item $(C,A)$ is observable;

\item for each $f \in \mathcal F$, $\mathcal G_{\mathcal O}$ is jointly connected by $\eta^f_{\mathcal O}$;

\item for each zero $z$ of $P_T(z)$, $z \neq 0$.
\end{enumerate}
\end{corollary}

\begin{corollary}
Given a MCN $\mathcal N$ and a faulty set $\mathcal F$, there exists a weight function $\bar W_{\mathcal O}$ such that $\mathcal N_f$ is detectable for each $f \in \mathcal F$, if and only if the following hold:
\begin{enumerate}
\item $(C,A)$ is detectable;

\item for each $f \in \mathcal F$, $\mathcal G_{\mathcal O}$ is jointly connected by $\eta^f_{\mathcal O}$.
\end{enumerate}

\end{corollary}

\section{Example: stabilization of a faulty MCN} \label{secSimulations}

Given the MCN $\mathcal N = ( \mathcal P, \mathcal G_{\mathcal R}, \eta_{\mathcal R}, \mathcal G_{\mathcal O}, \eta_{\mathcal O}, \Delta)$ defined in Example \ref{exRedundancy}, we show that the methodological results developed in this paper allow co-design of control algorithms and communication parameters for stabilizing a MCN. For clarity of presentation and without loss of generality, we assume that the set of failures configurations only takes into account the fault $f_3 = \{(v_3,v_7),(v_2,v_7)\}$, namely $\mathcal F = \{\varnothing, f_3\}$. As already illustrated in Example \ref{exRedundancy}, when the fault $f_3$ occurs in the links $(v_3,v_7)$ and $(v_2,v_7)$, then the zero $z_{f_3} = -2$ of $G^{f_3}_{\mathcal R}(z)$ coincides with the pole $p_2$ of $P_{T^c}(z)$, and the MCN $\mathcal N_{f_3}$ becomes uncontrollable. Applying Theorem \ref{thfaultyMCNControllability}, it is possible to replace the weight $W_{\mathcal R}(v_4,v_7)$ by $\bar W_{\mathcal R}(v_4,v_7) = W_{\mathcal R}(v_4,v_7)+0.1$, so that both $\mathcal N_\varnothing$ and $\mathcal N_{f_3}$ are controllable and observable.

Consider the system $\mathcal N$, that models the dynamical behavior of the network in case of failures. $\mathcal N$ can be modeled by a hybrid system as defined in \cite{BalluchiHSCC2002}, where the discrete-time dynamics switches from the non-faulty behavior $\mathcal N_\varnothing$ to the faulty behavior $\mathcal N_{f_3}$ and vice-versa. We define a minimum dwelling time $\tau$ for $\mathcal N$, such that the time duration between two consecutive changes of dynamics of $\mathcal N$ can not be smaller than $\tau$.

\begin{figure}[ht]
\begin{center}
\includegraphics[width=0.35\textwidth]{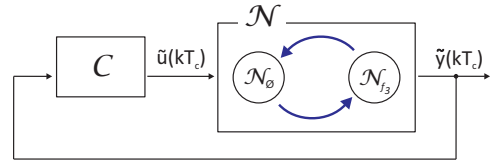}
\caption{Control scheme for the faulty MCN $\mathcal N$.} \label{figFaultyMCNblocks}
\end{center}
\end{figure}

In order to stabilize the closed loop system as depicted in Figure \ref{figFaultyMCNblocks}, we need to design a controller $\mathcal C$. However, it is not easy to guarantee the existence of a unique controller $\mathcal C$ by means of a discrete-time SISO LTI system, such that the switching dynamics of the closed loop $\mathcal N$ is stable. In order to design a stabilizing controller $\mathcal C$, we apply classical methodologies for eigenvalue placement to define a controller $\mathcal C_\varnothing$ for the non-faulty system $\mathcal N_\varnothing$, and a controller $\mathcal C_{f_3}$ for the faulty system $\mathcal N_{f_3}$. Since we designed the weight function $\bar W_{\mathcal R}$ so that $\mathcal N_\varnothing$ and $\mathcal N_{f_3}$ are controllable and observable, the existence of $\mathcal C_\varnothing$ and $\mathcal C_{f_3}$ is guaranteed. In order to decide for each time instant what control algorithm we have to apply to $\mathcal N$, we need a strategy to detect on-the-fly, and only using the input and output of $\mathcal N$, whether the current dynamics are non-faulty $\mathcal N_\varnothing$ or faulty $\mathcal N_{f_3}$.

\begin{figure}[b]
\begin{center}
\includegraphics[width=0.35\textwidth]{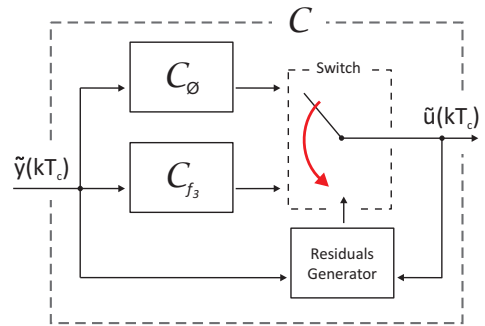}
\caption{Block diagram of the controller $\mathcal C$.} \label{figResidualGeneratorBlock}
\end{center}
\end{figure}

As illustrated in Figure \ref{figResidualGeneratorBlock} we used the hybrid observer techniques developed in \cite{BalluchiHSCC2002} and \cite{BalluchiIFAC2005}, which allow to construct a residual generator dynamical system able to detect the current dynamics of the system $\mathcal N$. The stability of the control scheme which exploits the residual generator is guaranteed to be stable, if the faulty and non-faulty dynamics switch with a dwelling time $\tau$ that is sufficiently smaller than the plant dynamics. This assumption is reasonable in our case study, since we consider link failures (e.g. due to battery discharge of a node), and thus characterized by rare occurrence.

\begin{figure}[t]
\begin{center}
\includegraphics[width=0.5\textwidth]{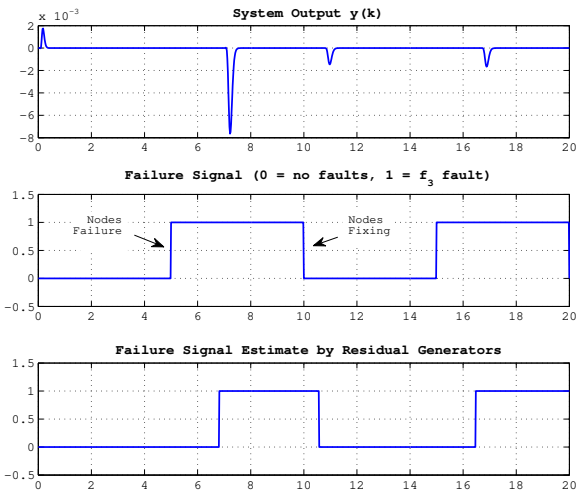}
\caption{Simulations results.} \label{figSimulations}
\end{center}
\end{figure}

We implemented the closed loop system in Simulink and performed the simulations in Figure \ref{figSimulations}, which show that when the dynamics of $\mathcal N$ switches, the controller detects the change of dynamics and applies the appropriate stabilizing control law.


\section{Conclusions} \label{secConclusions}

This work provides a novel methodology to design scheduling and routing of a communication network in order to preserve controllability and observability, for any set of failures configurations that at least preserve connectivity within the scheduling period between the controller and the plant, and vice-versa. In Section \ref{secSimulations} we showed that the configuration of scheduling and routing, together with the classical methodologies for eigenvalues placement, represent a novel co-design methodology of the network parameters and of the control algorithm for multi-hop control networks.

In future extensions of this paper, we aim to address the same problem for MIMO systems. Moreover, we will address the issue of introducing dynamical routing in our model, and performing an optimal choice of scheduling and weight functions. Another interesting problem to be addressed is guaranteeing the existence of a unique LTI controller of a MCN, that guarantees stability of the closed loop although the switching dynamics introduced by failures.


\bibliography{MCNControllability}
\end{document}